\documentclass[10pt]{amsart}
\usepackage{amssymb,amsmath,amsthm}
\usepackage{amssymb}
\usepackage{latexsym}
\usepackage{amscd}
\usepackage{color}
\input amssym.def
\input amssym
\newtheorem{theorem}{Theorem}[section]
\newtheorem{lemma}[theorem]{Lemma}
\newtheorem{definition}[theorem]{Definition}

\newtheorem{proposition}[theorem]{Proposition}

\newtheorem{remark}{Remark}

\newcommand{\W}{\frac{W'_V(\tau)}{W_V(\tau)}}
\newcommand{\WC}{\frac{\mathcal{W}'_V(\tau)}{\mathcal{W}_V(\tau)}}

\newcommand{\halmos}{\rule{1ex}{1.4ex}}

\newcommand{\bea}{\begin{eqnarray}}
\newcommand{\epf}{\hspace*{\fill}\mbox{$\halmos$}}
\newcommand{\eea}{\end{eqnarray}}
\newcommand{\nn}{\nonumber \\}
\newcommand{\be}{\begin{equation}}

\newcommand{\ee}{\end{equation}}

\usepackage{amssymb}

\begin{document}
\title[Modular identities and supersingular $j$-invariants]{Modular invariance, modular identities and supersingular $j$-invariants}
\author{Antun  Milas}
\thanks{}
\begin{abstract}
To every $k$-dimensional modular invariant vector space we
associate a modular form on $SL(2,\mathbb{Z})$ of weight $2k$. We
explore number theoretic properties of this form and find a
sufficient condition for its vanishing which yields modular
identities (e.g., Ramanujan-Watson's modular identities).
Furthermore, we focus on a family of modular invariant spaces
coming from suitable two-dimensional spaces via the symmetric
power construction. In particular, we consider a two-dimensional
space spanned by graded dimensions of certain level one modules
for the affine Kac-Moody Lie algebra of type $D_4^{(1)}$. In this
case, the reduction modulo prime $p=2k+3 \geq 5$ of the modular
form associated to the $k$-th symmetric power classifies
supersingular elliptic curves in characteristic $p$. This
construction also gives a new interpretation of certain modular
forms studied by Kaneko and Zagier.
\end{abstract}

\address{Department of Mathematics and Statistics,
SUNY Albany, Albany, New York 12222}
\email{amilas@math.albany.edu}

\subjclass[2000]{11F03, 11F30, 17B67, 17B69} \maketitle

\renewcommand{\theequation}{\thesection.\arabic{equation}}
\renewcommand{\thetheorem}{\thesection.\arabic{theorem}}
\setcounter{equation}{0} \setcounter{theorem}{0}
\setcounter{equation}{0} \setcounter{theorem}{0}
\setcounter{section}{0}

\section{Introduction and notation}
An especially interesting feature of every rational vertex
operator algebra is the modular invariance of graded dimensions
(see \cite{Zh} for a precise statement).
What distinguishes modular invariant spaces coming from
representations of vertex operator algebras is the fact that these
spaces are equipped with special spanning sets indexed by
irreducible modules of the algebra, and are subject to the
Verlinde formula (cf. \cite{Hu}). Moreover, every irreducible
graded dimension (or simply, {\em character}) admits a
$q$-expansion of the form
$$q^{\bar{h}} \sum_{n=0}^\infty a_n q^n,$$
where $\bar{h} \in \mathbb{Q}$ and $a_n \in \mathbb{Z}_{\geq 0}.$

In \cite{M1} we showed that the internal structure of certain
vertex operator algebras can be conveniently used to prove some
modular identities without much use of the theory of modular
forms. The key ingredient in our approach was played by certain
Wronskian determinants which are intimately related to ordinary
differential equations with coefficients being holomorphic modular
forms. Additionally, these differential equations are closely
related to certain finiteness condition on the vertex operator
algebra in question. In a joint work with Mortenson and Ono
\cite{MMO} we studied differential equations associated to
$(2,2k+1)$ Andrews-Gordon series and observed that a suitably
normalized constant coefficient in these ODEs (expressible as a
quotient of two Wronskians), when restricted modulo prime $p=2k+1$
is essentially the locus of supersingular $j$-invariants in
characteristic $p$. It is an open problem to find an alternative
description of the modular forms considered in \cite{MMO}.

The aim of this note is to build a framework for studying modular
forms expressed as a quotient of two Wronskians as in \cite{MMO}.
In addition we present a very elegant way of constructing
supersingular polynomials by using symmetric products. Thus we are
able to explicitly determine our modular forms and to relate them
to known constructions in the literature (cf. \cite{KZ}).

First we show how to construct a modular form of weight $2k$ from
a $k$-dimensional modular invariant space. Let $V$ be {\em modular
invariant} vector space with a basis $f_1(\tau),...,f_k(\tau)$
\footnote{We refer to Section 2 for precise conditions on
$f_i(\tau)$.},
i.e., for every $i$ and $\gamma=\left[ \begin{array}{cc} a & b \\
c & d \end{array} \right] \in SL(2, \mathbb{Z})$, there exist
constants $\gamma_{i,j}$ such that
$$f_i \left(\frac{a \tau+b}{c \tau+d} \right)=\sum_{j=1}^k \gamma_{i,j}
f_j(\tau).$$ Let us denote the Wronskian determinant of
$f_1(\tau)$,...$f_k(\tau)$ by
$$W_V=W_{(q \frac{d}{dq})}(f_1,...,f_k),$$
where we use the Ramanujan's derivative $\left(q \frac{d}{dq}
\right)$. This is an automorphic form on $SL(2,\mathbb{Z})$ (more
precisely, a modular form with a character) of weight $k(k-1)$ and
its properties have been recorded in the literature (cf.
\cite{AM}, \cite{M1}, \cite{M2}, \cite{MMS}, \cite{Mu1}). It is
not hard to see that the sixth power of $W_V$ is a modular form.
In addition, the Wronskian of the {derivatives} of
$f_1(\tau),...,f_k(\tau)$,
$$W'_V=W_{(q \frac{d}{dq})}(f'_1,...,f'_k),$$ is also an
automorphic form of weight $k(k+1)$ with the same character as
$W_V$ (cf. \cite{M2}). Unlike $W_V$ and $W'_V$, \be \label{qq} \W
\ee is independent of the particular choice of the basis of $V$.
This quotient, which is the main object of our study, is a
(meromorphic) modular form for $SL(2,\mathbb{Z})$ of weight $2k$
(possibly zero) \cite{M2}, \cite{MMO}.
We also denote by $\mathcal{W}_V$ (resp. $\mathcal{W}'_V$) the
normalization of $W_V$ (resp. $W'_V$) (if nonzero) in which in the
$q$-expansion the leading coefficient is one. Clearly,
$\mathcal{W}_V$ and $\mathcal{W}_V$ do not depend on the basis
chosen.

Alternatively, we can think of $\frac{W'_V}{{W}}$ as follows.
There is a unique linear differential operator $\mathcal{D}_V$ of
order $k$,
$$\mathcal{D}_V=\sum_{i=0}^k P_{i,V}(q) \left(q \frac{d}{dq} \right)^i,$$
which satisfies $\mathcal{D}_V y=0$  for every $y \in V$, and
$P_{k,V}=1.$ Under these conditions
$$P_{0,V}(\tau)=(-1)^k \W,$$
so $\frac{{W'_V}}{{W}}$ is , up to a sign, just the evaluation
$\mathcal{D}_V(1)$.

This paper is organized as follows.
In Section 2 we obtain a sufficient condition for the vanishing of
$\W$ (see Theorem \ref{one}). Then in Section 3 we focus on
modular invariant spaces obtained via the { symmetric power}
construction (cf. Theorem \ref{2sym}). In Section 4 we apply
results from sections 2 and 3 to prove some modular identities,
such as the Ramanujan-Watson's modular identities for the
Rogers-Ramanujan's continued fraction (cf. Theorem \ref{onet}). In
Section 5 we derive a recursion formula (cf. Lemma \ref{prelast})
which can be used to give another proof of Theorem \ref{onet}. We
derive the same recursion in Section 6 in the framework of vertex
operator algebras (this part can be skipped without any loss of
continuity). In Section 7 we gather some results about
supersingular elliptic curves and modular forms. Finally, in
Section 8 we focus on supersingular congruences for modular forms
obtained from $V={\rm Sym}^m(U)$, where $U$ is spanned by the
graded dimensions of level one modules for $D_4^{(1)}$ (see also
\cite{MMO} for a related work). Our main result, Theorem
\ref{two}, gives a nice expression for $\W$ as a coefficient of a
certain generating series studied in \cite{KZ}.

Throughout the paper the Eisenstein series will be denoted by \be
\label{g2k} G_{2k}(\tau)=\frac{-B_{2k}}{2k!}+\frac{2}{(2k-1)!}
\sum_{n =1}^\infty \frac{n^{2k-1}q^n}{1-q^n}, \ k \geq 1. \ee We
will also use normalized Eisenstein series \be \label{e2k}
E_{2k}(\tau)=(\frac{-B_{2k}}{2k!})^{-1}G_{2k}(\tau). \ee

As usual, the Dedekind $\eta$-function and the discriminant are
defined as \bea \eta(\tau)&=&q^{1/24} \prod_{n=1}^\infty (1-q^n)
\nn \Delta(\tau)&=&E_4(\tau)^3-E_6(\tau)^2. \eea Then the
$j$-function is defined as
$$j(\tau)=\frac{1728 E_4(\tau)^3}{\Delta(\tau)}.$$
A holomorphic modular form for $SL(2,\mathbb{Z})$ is assumed to be
holomorphic in $\mathbb{H}$ (the upper half-plane) with a possible
pole at the infinity. The order of vanishing at the infinity of
$f$ will be denoted by ${\rm ord}_{i \infty}(f)$. The graded ring
of holomorphic modular forms including at the infinity will be
denoted by $M=\mathbb{C}[E_4,E_6]$. Its graded components will be
denoted by $M_k$. Every $f(\tau) \in M_{k}$ can be written
uniquely as \be \label{rep} f(\tau)=\Delta^t(\tau)
E_4^\delta(\tau) E_6^\epsilon(\tau) \tilde{F}(f,j(\tau)), \ee
where $\tilde{F}(f,j)$ is a polynomial of degree $\leq t$ and
$$k=12t+4 \delta+6 \epsilon,$$
where $0 \leq \delta \leq 2$ and $0 \leq \epsilon \leq 1$.

\section{On the vanishing of $\W$}

\noindent The goal of this section is to obtain a sufficient
condition for the vanishing of $\W$. All our modular invariant
spaces are assumed to have a basis $\{f_1(\tau),...,f_k(\tau) \}$
of holomorphic functions in $\mathbb{H}$ with the $q$-expansion of
the form \be \label{fi} f_i(\tau)=q^{h_i} \sum_{n=0}^\infty
a_n^{(i)} q^n, \ee where $h_i \in \mathbb{Q}$. The rationale for
this assumption rests on the general form of irreducible
characters of rational vertex operator algebras \cite{DLM} (see
also \cite{AM}). In fact, in all our applications $a_n^{(i)}$ are
nonnegative integers. We start with an auxiliary result:
\begin{lemma} \label{vanlemma} Let $\{ f_1(\tau),...,f_k(\tau) \}$ be a basis of $V$ with
$q$-expansions as in (\ref{fi}). Then we can find a (new) basis of
$V$ \be \label{qexp} \bar{f}_i(\tau)=q^{\bar{h}_i}
\sum_{n=0}^\infty \bar{a}_n^{(i)} q^n, \ \ \bar{a}_0^{(i)} \neq 0,
\ i=1,...,k, \ee where \be \label{exp} \bar{h}_1<\bar{h}_2< \cdots
< \bar{h}_k, \ee so that
$${\rm ord}_{i \infty} W_V(\tau)=\sum_{i=1}^k \bar{h}_i.$$
The numbers $\bar{h}_i$ are uniquely determined.
\end{lemma}
\noindent {\em Proof:} The uniqueness of $\bar{h}_i$ follows
easily by the induction on $k$. Clearly, $W_{(q
\frac{d}{dq})}(\bar{f}_1,...,\bar{f}_k)$ is a nonzero multiple of
$W_{(q \frac{d}{dq})}({f}_1,...,{f}_k)$. Now, the leading
coefficient in the $q$-expansion of $W_V$ is (up to a sign) the
Vandermonde determinant $V(\bar{h}_1,...,\bar{h}_k)=\prod_{i <j}
(\bar{h}_i-\bar{h}_j) \neq 0$, and the leading power of $q$ is
$\sum_{i=1}^k \bar{h}_i$. \epf

The vanishing of $W'_V$ simply means that there is a linear
relation \be \label{C} \sum_{i=1}^k \lambda_i f_i(\tau)=C \neq 0.
\ee The following result gives a sufficient condition for the
vanishing of $W'_V$.
\begin{theorem} \label{one}
Let $\{ f_1(\tau),...,f_k(\tau) \}$ be a basis of a modular
invariant space $V$ satisfying (\ref{qexp}) and (\ref{exp}) such
that
\begin{itemize}
\item[(i)] $\frac{{W}'_V(\tau)}{W_V(\tau)}$ is holomorphic (i.e.,
$W_V(\tau)$ is nonvanishing in $\mathbb{H}$). \item[(ii)] There
exists $r \geq \lfloor \frac{k}{6} \rfloor$ and $
f_{i_0}(\tau),...,f_{i_r}(\tau)$ with \be {\rm ord}_{i
\infty}f_{i_j}(\tau)=j, \ \ {\rm for} \ j=0,...,r. \ee
\end{itemize}
Then $\W$ is identically zero. If in addition $
f_{i_0}(\tau)$,...,$f_{i_r}(\tau)$ are the only $f_i$ with
positive integer powers of $q$, then there exist constants
$\lambda_i$, $C$, such that
$$\sum_{j=0}^r \lambda_i f_{i_j}(\tau)=C \neq 0.$$
\end{theorem}
\noindent {\em Proof:} Clearly, $\bar{h}_{i_j}=j$ for $j=0,...,r$.
Because of ${\rm ord}_{i \infty} f'_{i_0}(\tau) \geq 1$ and ${\rm
ord}_{i \infty} f'_{i_j}(\tau)=j$ for $j \geq 1$, we can find
constants $\lambda_j$ such that \be {\rm ord}_{i \infty}
(\sum_{j=0}^r \lambda_j f'_{i_j}(\tau)) \geq r+1. \ee We claim
that $\W$ is zero. Suppose that $\W \neq 0$. By Lemma
\ref{vanlemma}, we have ${\rm ord}_{i \infty}
W_V(\tau)=\sum_{i=1}^k \bar{h}_i$ and ${\rm ord}_{i \infty}
W'_V(\tau) \geq r+1+ \sum_{i=1}^k \bar{h}_i$ (keep in mind that
$W'_V$ is just the Wronskian of $f'_1(\tau),...,f'_k(\tau)$).
Thus,
$${\rm ord}_{i \infty} \frac{W'_V(\tau)}{W_V(\tau)} \geq r+1
> \lfloor \frac{k}{6} \rfloor.$$ It is known that the order of
vanishing at the infinity of a nonzero holomorphic modular form of
weight $2k$ is at most $\lfloor \frac{k}{6} \rfloor$. The first
claim holds.

Suppose now that $f_{i_0}(\tau)$,...,$f_{i_r}(\tau)$ are the only
$f_i$ with positive integer powers of $q$. Because of (\ref{qexp})
$q$-powers of $f_i$ are integral if and only if $\bar{h}_i$ is an
integer. If $\sum_{j=0}^r \lambda_j f'_{i_j}(\tau)$ is nonzero,
then ${\rm ord}_{i \infty} \sum_{j=0}^r \lambda_j f'_{i_j}(\tau)$
is finite and therefore ${\rm ord}_{i \infty} W'_V$ is also finite
and $W'_V$ is nonzero. We have a contradiction. \epf


\section{Wronskians and symmetric powers}

\begin{definition} {\em Let $U$ be a modular invariant space and $m$ a positive integer.
The modular invariant space spanned by
$$ \{ f_1 \cdots f_m : \ f_i \in U \},$$
is called the $m$-th symmetric power \footnote{This terminology
will be explained later in Section 5.} of $U$ and is denoted by
${\rm Sym}^m(U)$.}
\end{definition}
For ${\rm dim}(U) \geq 3$ it is a nontrivial task to find even the
dimension of $V={\rm Sym}^m(U)$, let alone to extract any
information regarding $\W$. However, if ${\rm dim}(U)=2$, the
situation is much better and we have the following result.
\begin{theorem} \label{2sym}
Let $U$ be a two-dimensional modular invariant space, then $V={\rm
Sym}^{m}(U)$ is $(m+1)$-dimensional and
$$\mathcal{W}_V(\tau)=\mathcal{W}_U(\tau)^{\frac{m(m+1)}{2}}.$$
If in addition $W_U(\tau)$ is nonvanishing, then $\W$ is a
holomorphic modular form of weight $2m+2$ and \be \label{eta}
\mathcal{W}_V(\tau)=\eta(\tau)^{2m(m+1)}. \ee
\end{theorem}
\noindent {\em Proof:} Let $f_1$ and $f_2$ form a basis of $U$.
Then the set \be \label{sym2} \{ f_1^i f_2^{m-i} : \ i=0,...,m \},
\ee is linearly independent (otherwise $f_1/f_2$ would be a
constant), so it gives a basis of $V$. Now, by using basic
properties of the Wronskian, we have \bea W_V(\tau)&=&W_{(q
\frac{d}{dq})}(f_1^m,f_1^{m-1}f_2,...,f_1f_2^{m-1},f_2^m) \nn &=&
(f_1^m)^{m+1} W_{(q \frac{d}{dq})}(1,(f_2/f_1),...,(f_2/f_1)^m)=
f_1^{m(m+1)} W_{(q \frac{d}{dq})}((f_2/f_1)',...,((f_2/f_1)^m)')
\nn &=& f_1^{m(m+1)} ((f_2/f_1)')^m W_{(q \frac{d}{dq})}(
1,(f_2/f_1),...,(m-1)(f_2/f_1)^{m-1}) \nn &=& m! f_1^{m(m+1)}
((f_2/f_1)')^m W_{(q
\frac{d}{dq})}(1,(f_2/f_1),...,(f_2/f_1)^{m-1}) = \left(
\prod_{k=1}^m k! \right) f_1^{m(m+1)} ((f_2/f_1)')^{m(m+1)/2} \nn
&=&\left( \prod_{k=1}^m k! \right)
(f_2'f_1-f_1'f_2)^{m(m+1)/2}=\left(\prod_{k=1}^m k!
\right)W_U(\tau)^{m(m+1)/2}. \nonumber \eea If $W_U$ is
nonvanishing then $\mathcal{W}_U(\tau)=\eta(\tau)^4$ (cf.
\cite{M2}) and (\ref{eta}) follows. \epf


\section{Ramanujan-Watson's modular identities} \noindent

\noindent The Rogers-Ramanujan continued fraction \cite{A} is
defined as
$$R(q):=\frac{q^{1/5}}{1+ \displaystyle{\frac{q}{\displaystyle{1+\frac{q}{1 +
\cdots}}}}}.$$ In one of his notebooks Ramanujan stated that
$R(e^{-\pi \sqrt{r}})$ can be exactly found for every positive
rational number $r$. The main identities that support Ramanujan's
claim are the Rogers-Ramanujan identities (cf. \cite{A})  and a
pair of modular identities recorded by Ramanujan \cite{Ra}:
\begin{theorem} \label{onet}
We have \be \label{wa}
\frac{1}{R(q)}-1-{R(q)}=\frac{\eta(\tau/5)}{\eta(5 \tau)}, \ee \be
\label{wat2}
\frac{1}{R(q)^5}-11-R(q)^5=\left(\frac{\eta(\tau)}{\eta(5
\tau)}\right)^6. \ee
\end{theorem}
The first proof of Theorem \ref{onet} was obtained by Watson
\cite{Wa}. There are other proofs in the literature that use
methods similar to those available to Ramanujan (see \cite{B} for
another proof and \cite{BHSS} for a discussion on this subject).
More analytic proofs use nontrivial facts such as explicit forms
of Hauptmodules for certain modular curves (see \cite{D} for a
nice review). We will prove (\ref{wat2}) by using Theorem
\ref{one}.

Firstly, we will need the following well-known fact (we refer the
reader to \cite{CCM} for a discussion in the context of the
two-dimensional conformal field theory):
\begin{lemma} \label{twot}
Let $U$ be the vector space spanned by \be \label{rr1}
ch_1(\tau):=q^{11/60} \prod_{n=0}^\infty
\frac{1}{(1-q^{5n+2})(1-q^{5n+3})},\ee \be \label{rr2}
ch_2(\tau):=q^{-1/60} \prod_{n=0}^\infty
\frac{1}{(1-q^{5n+1})(1-q^{5n+4})}. \ee The modular transformation
$$\tau \mapsto \frac{-1}{\tau},$$
induces an endomorphism of $U$, which in the basis $\{ ch_1, ch_2
\}$ is represented by the matrix
$$S=\frac{2}{\sqrt{5}} \left[\begin{array}{cc} -sin \left(\frac{2
\pi}{5} \right) & sin \left(\frac{4 \pi}{5} \right) \\ sin
\left(\frac{4 \pi}{5} \right) & sin \left(\frac{2 \pi}{5} \right)
\end{array} \right].$$
\end{lemma}
\noindent {\em Proof:} By using Jacobi Triple Product Identity we
first rewrite (\ref{rr1})-(\ref{rr2}) as quotients of two theta
constants \bea \label{last1} ch_1(\tau)&=&\frac{q^{9/40} \sum_{n
\in \mathbb{Z}} (-1)^n q^{\frac{5n^2+3n}{2}}}{\eta(\tau)}, \\
\label{last2} ch_2(\tau)&=&\frac{q^{1/40} \sum_{n \in \mathbb{Z}}
(-1)^n q^{\frac{5n^2+n}{2}}}{\eta(\tau)}. \eea  Now, apply the
formula \be \label{ded} \eta(-1/\tau)=\sqrt{-i \tau } \eta(\tau).
\ee and the modular transformation formulas for the two theta
constants in the numerators of (\ref{last1})-(\ref{last2}), under
$\tau \mapsto \frac{-1}{\tau}$. For an explicit computations in this
case see, for instance, \cite{D}.
\epf

\noindent {\em Proof of (\ref{wa}):}
Observe first that
$$ch_1 (\tau) \cdot ch_2 (\tau)=\frac{\eta(5 \tau)}{\eta(\tau)}.$$
From (\ref{ded}) we have \be \label{eq1} ch_1(-1/\tau)
ch_2(-1/\tau)=\frac{\sqrt{-i \tau/5} \eta(\tau/5)}{\sqrt{-i \tau}
\eta(\tau)}=\frac{1}{\sqrt{5}} \frac{\eta(\tau/5)}{\eta(\tau)}.
\ee On the other hand because of the lemma and a few trigonometric
identities for $sin \left(\frac{2 \pi}{5} \right)$ and $sin
\left(\frac{4 \pi}{5}\right)$ \be \label{eq2} ch_1(-1/\tau)
ch_2(-1/\tau) =\frac{1}{\sqrt{5}} \left( -ch_1(\tau)
ch_2(\tau)-ch_1(\tau)^2+ch_2(\tau)^2 \right). \ee Now, after we
equate the right-hand sides of (\ref{eq1}) and (\ref{eq2}), cancel
the factor $\frac{1}{\sqrt{5}}$ and multiply both sides by
$$\frac{1}{ch_1(\tau) ch_2(\tau)}=\frac{\eta(\tau)}{\eta(5
\tau)},$$ we get \be \label{ww}
-1-\frac{ch_1(\tau)}{ch_2(\tau)}+\frac{ch_2(\tau)}{ch_1(\tau)}=\frac{\eta(\tau/5)}{\eta(5
\tau)}. \ee Finally, we recall the Rogers-Ramanujan identities
\cite{A}:
$$R(q)=q^{1/5} \frac{\prod_{n=0}^\infty
(1-q^{5n+1})(1-q^{5n+4})}{\prod_{n=0}^\infty
(1-q^{5n+2})(1-q^{5n+3})}$$ and observe that
$R(\tau)=\frac{ch_1(\tau)}{ch_2(\tau)}$. \epf

%
%
\noindent {\em Proof of (\ref{wat2}):} We will prove the following
equivalent statement:
$$\frac{ch_2(\tau)^5}{ch_1(\tau)^5}-11-\frac{ch_1(\tau)^5}{ch_2(\tau)^5}=\left(\frac{1}{ch_1(\tau)
ch_2(\tau)}\right)^6,$$ which can be rewritten as \be
\label{wa2-ch} ch_2(\tau)^{11} ch_1(\tau)-11 ch_1(\tau)^6
ch_2(\tau)^6-ch_1(\tau)^{11} ch_2(\tau)=1. \ee Consider the
$12$-th symmetric power of $U$ with a basis
$$ \{ ch_1^i(\tau) ch_2^{12-i}(\tau), \ \ 0 \leq i \leq 12 \}.$$
If we let $ch_1^i(\tau) ch_2^{12-j}(\tau)=q^{\bar{h}_j}+\cdots,$
then the exponents $\bar{h}_j$ satisfy $\bar{h}_0<\bar{h}_1<
\cdots < \bar{h}_{12}$. The crucial observation here is that
$$ch_1^i (\tau) ch_2^{12-i}(\tau) \in \mathbb{Q}[[q]],$$
if and only if $i=1$, $i=6$ or $i=11$. For all other $i$ the
powers are nonintegral. More precisely, \bea ch_1(\tau)
ch_2^{11}(\tau)&=&1+11 q+ 67q^2 \cdots, \nn ch_1^6(\tau)
ch_2^6(\tau)&=&q+6q^2+\cdots,\nn ch_1^{11}(\tau)
ch_2(\tau)&=&q^2+\cdots. \eea

Observe that \be {\rm ord}_{i \infty} (\label{11}
({ch}_{1}^{}(\tau)ch_2^{11}(\tau))'-11
(ch_1^6(\tau)ch_2^6(\tau))'-(ch_1^{11}(\tau)ch_2(\tau))') \geq 3.
\ee Now, we are ready to apply Theorem \ref{one}. Here $k=13$,
$r=2$, and $\mathcal{W}_{{\rm
Sym}^{12}(U)}(\tau)=\Delta(\tau)^{13}$ is nonvanishing (cf.
\cite{M1} and Theorem \ref{2sym}). Now, (\ref{11}) and Theorem
\ref{one} imply \be \label{wa3-ch} ch_2(\tau)^{11} ch_1(\tau)-11
ch_1(\tau)^6 ch_2(\tau)^6-ch_1(\tau)^{11} ch_2(\tau)=C \neq 0. \ee
The constant $C$ is clearly $1$.
\epf

\begin{remark} {\em It is possible to prove (\ref{wat2}) without referring to Theorem \ref{one} and \cite{M1}.
Notice that (\ref{11}) implies that ${\rm ord}_{i \infty} (W'_V)
\geq 16.$ But there is no modular form of weight $13 \cdot 14=182$
with this behavior at the cusp.}
\end{remark}

Our Theorem \ref{one} can be applied in a variety of situations as
long as the degree of the symmetric power is not too big. Here we
apply our method in the case of the vector space spanned by \bea
f_1(\tau)&=&\frac{\sum_{n \in \mathbb{Z}} q^{n^2}}{\eta(\tau)},
\nn f_2(\tau)&=&\frac{\sum_{n \in \mathbb{Z}}
q^{(n+1/2)^2}}{\eta(\tau)}. \nonumber \eea Those familiar with the
theory of affine Kac-Moody Lie algebras will recognize these
series as modified graded dimensions of two distinguished
irreducible representations of the affine Lie algebra of type
$A_1^{(1)}$ \cite{K}. It is not hard to see that the vector space
spanned by $f_1(\tau)$ and $f_2(\tau)$ is modular invariant
\cite{CCM}, \cite{K}. Then we have the following analogue of
Theorem \ref{onet}:
\begin{proposition} \label{weber}
We have \be \label{ex1}
2\frac{f_1(\tau)}{f_2(\tau)}-2\frac{f_2(\tau)}{f_1(\tau)}=\left(\frac{\eta(\tau/2)}{\eta(2
\tau)}\right)^4,\ee and \be \label{ex2} f_1(\tau)^5
f_2(\tau)-f_2(\tau)^5 f_1(\tau)=2. \ee The identity (\ref{ex2}) is
equivalent to the following classical identity for Weber modular
functions: \be \label{ex3} \prod_{n=1}^\infty (1+q^{2n-1})^8-16q
\prod_{n=1}^\infty (1+q^{2n})^8=\prod_{n=1}^\infty (1-q^{2n-1})^8.
\ee
\end{proposition}
\noindent {\em Proof:} Firstly, we apply the Jacobi Triple Product
Identity \cite{A} so that \be \label{ex4}
f_1(\tau)=\frac{\prod_{n=1}^\infty
(1-q^{2n})(1+q^{2n-1})^2}{\eta(\tau)}, \ee \be \label{ex5}
f_2(\tau)=\frac{2 q^{1/4} \prod_{n=1}^\infty
(1-q^{2n})(1+q^{2n})^2}{\eta(\tau)}. \ee For (\ref{ex1}), notice
that
$$f_1(\tau) f_2(\tau)=2 \left(\frac{\eta(2
\tau)}{\eta(\tau)}\right)^4.$$ Now apply $\tau \mapsto
\frac{-1}{\tau}$ and proceed as in the proof of (\ref{wa}).

Similarly, (\ref{ex2}) follows from analysis of  $f_1^i(\tau)
f_2^{6-i}(\tau)$, by following the steps as in the proof of
(\ref{wat2}). The identity (\ref{ex3}) is now a consequence of
(\ref{ex2}), (\ref{ex4}) and (\ref{ex5}). \epf

\section{Vanishing of $\WC$: ODE approach}

We have seen that $\W=0$ can be deduced by a careful analysis of
the order of vanishing of $\W$ (or $W'_V$) at the infinity. In
this section we deduce a related result by using elementary theory
of ordinary differential equations.

To every second order homogenous ODE of the form \be \label{ode2}
\left(q \frac{d}{dq}\right)^2 y+P(q)\left(q \frac{d}{dq}
\right)y+Q(q) y=0, \ee we associate its $m$-th symmetric power ODE
which is by the definition a homogeneous ODE of minimal order with
a fundamental system of solutions
$$\{ f^i g^{m-i} \ i=0,...,m \},$$ where $\{ f, g  \}$ is a fundamental system of
solutions of (\ref{ode2}). See \cite{S} for more about symmetric
powers of ODEs in general. Now, unlike symmetric powers for
equations of the order three, the $m$-th symmetric power of
(\ref{ode2}) is always of order $m+1$ (cf. Section 3) and is given
by
$$(-1)^{m+1}
\frac{{W}_{(q \frac{d}{dq})}(y,
f^m,f^{m-1}g,...,fg^{m-1},g^m)}{{W}_{(q
\frac{d}{dq})}(f^m,f^{m-1}g,...,fg^{m-1},g^m)}=0.$$ By expanding
the determinant in the numerator we obtain \be \label{odep} \left(
q \frac{d}{dq} \right)^{m+1} y + \sum_{i=0}^{m} Q_{m,i}(q) \left(q
\frac{d}{dq}\right)^i y=0. \ee Clearly, the "constant" coefficient
$Q_{m,0}(q)$ is equal to \be \label{m1} (-1)^{m+1} \frac{{W}_{(q
\frac{d}{dq})}((f^m)',(f^{m-1}g)',...,(fg^{m-1})',(g^m)')}
{{W}_{(q \frac{d}{dq})}(f^m,f^{m-1}g,...,f g^{m-1},g^m)}. \ee We
will show that it is possible to compute $Q_{m,0}(q)$ via a
certain recursion formula. Let \be \label{thetafirst}
\Theta_h:=\left(q \frac{d}{dq} \right)+ h G_2(q). \ee
Then $\Theta_h : L_k \longrightarrow L_{k+2}$, where $L_k$ stands
for any modular invariant space of weight $k$ for $SL(2,
\mathbb{Z})$, in particular, the vector space of holomorphic
modular forms of weight $h$. We will use notation
$$\Theta^k:=\Theta_{2k} \circ \cdots \circ \Theta_2 \circ
\Theta_0.$$ From now on we will focus on the following ODE: \be
\label{stuff} \Theta^2 y+Q(q) y=0, \ee where $Q(q)$ is a
(meromorphic) modular form of weight 4.
\begin{lemma} \label{prelast}
Fix $m  \geq 2$. Let \bea \label{recursion} R_1&=& m Q, \nn
R_2&=&m \Theta Q, \nn R_{i+1}&=&\Theta R_i+(i+1)(m-i) Q R_{i-1}, \
\ i=2,...,m-1. \eea Then
$$Q_{m,0}=R_m.$$
\end{lemma}
\noindent {\em Proof:}
Fix a positive integer $m$. The $m$-th symmetric power of
(\ref{stuff}) is given by
$$D_{m+1}y=0,$$ where $D_{m+1}$ is obtained recursively from \bea \label{recdiff}
D_0&=&1, \nn D_1&=& \Theta \nn \label{recursion2} D_{i+1}&=&\Theta
D_i+i(m-i+1) Q(q) D_{i-1}, \ \ 0 < i \leq m. \eea This recursion
formula can be proven by induction and seems to be known in the
literature (see for instance Theorem 5.9 in \cite{Do}). For
example, the second symmetric power ($m=2$) of (\ref{stuff}) is
given by
$$ D_3y=0, \ \ {\rm where} \ \ D_3=\Theta^3+4 Q \Theta+2 \Theta (Q).$$
Since every differential operator $D_i$ (which depends on $m$)
admits an expansion
$$D_i=\sum_{j=0}^{i} R_{j,i}(q) \Theta^i,$$
if we let now
$$R_{j-1}=R_{j,0}, \ \ j \geq 2,$$
then from the formula (\ref{recursion2}) we have \bea
R_1&=&R_{2,0}=mQ, \nn R_2&=&R_{3,0}=m \Theta Q, \nn
R_{i+1}&=&R_{i,0}=\Theta R_{i}+(i+1)(m-i) Q R_{i-1}, \ \ 2 \leq i
\leq m-1. \nonumber \eea \epf

\noindent Now, we specialize everything to an ODE of type \be
\label{ODEl} \Theta^2 y+\lambda G_4(\tau)y=0, \ \ \lambda \in
\mathbb{C}. \ee

\begin{lemma} \label{lemmalast}
Let $m=12$ and $Q=\lambda G_4(\tau)$. Then \be R_{12}=0, \ \
\mbox{if and only if} \ \ \lambda \in \{ -\frac{11}{5},
-\frac{25}{4}, -15,-40, 0 \}. \ee
\end{lemma}
\noindent {\em Proof:} Follows after some computation by using
Lemma \ref{prelast} and the formulas \bea \Theta G_4&=& 14 G_6 ,
\nn \Theta G_6&=& \frac{60 G^2_4}{7}, \eea known to Ramanujan.
\epf

The following proposition is from \cite{MMS} (it was also proven
in  \cite{M1}):
\begin{proposition} \label{milas}
The series ${ch}_1(\tau)$ and ${ch}_2(\tau)$ form a fundamental
system of solutions of \be \label{milas1} \Theta^2 y-\frac{11}{5}
G_4(\tau)y=0. \ee
\end{proposition}

\noindent {\em Proof of (\ref{wat2}):} The Proposition \ref{milas}
and the vanishing of $Q_{12,0}(q)$ in Lemma \ref{lemmalast} for
$\lambda=-\frac{11}{5}$ implies the vanishing of $\W$.
The proof now follows. \epf

\section{The recursion (\ref{recursion}) via vertex operator
algebras}

In \cite{M1} we obtained a representation theoretic proof of a
pair of Ramanujan's identities based on an internal structure of
certain irreducible representations of the Virasoro algebra. The
same framework can be used to prove the formula (\ref{wat2}).

In this section we will use the notation from \cite{M1}. Let $U$
be as in Section 4. Notice that ${\rm Sym}^{12}(U)$ is just the
vector space spanned by graded dimensions of irreducible modules
of the tensor product vertex operator algebra
$L(-22/5,0)^{\otimes^{12}}$ \cite{FHL}, where $L(-22/5,0)$ is the
vertex operator algebra associated to $\mathcal{M}(2,5)$ Virasoro
minimal models \cite{FF} \cite{M1}.
Let $L(n)$ and $L[n]$ be two sets of generators of the Virasoro
algebra as in \cite{Zh}, \cite{M1}. We showed in \cite{M1} that
\be \label{l20} {\rm tr}|_{W} o(L[-2]^2{\bf 1}) q^{L(0)+11/60}=0,
\ee for every $L(-22/5,0)$--module $W$.
For $ 0 \leq i \leq 12$, let
$$v_{i}=i! S(L[-2]{\bf 1} \otimes \cdots \otimes L[-2]{\bf 1} \otimes {\bf 1}
\otimes \cdots \otimes {\bf 1}) \in L(-22/5,0)^{\otimes^{12}},$$
where in the first $i$ tensor slots we have the vector $L[-2]{\bf
1}$, and on the remaining $(12-i)$ tensor slots the vector ${\bf
1}$, and $S$ denotes the symmetrization (e.g., $S(L[-2] \otimes
L[-2] \otimes {\bf 1})=L[-2] \otimes L[-2] \otimes {\bf 1}+L[-2]
\otimes {\bf 1} \otimes L[-2] +{\bf 1} \otimes L[-2] \otimes
L[-2]$).
It is known (see for instance \cite{Zh}, \cite{M1}), that for
every vertex operator algebra $V$, a $V$-module M, and a
homogeneous vector $w$ the following identity holds $${\rm tr}|_M
o(L[-2]w)q^{L(0)-c/24}$$
$$=\left(\left( q \frac{d}{dq} \right)+{\rm deg}(w)G_2(\tau)
\right) {\rm tr}|_M o(w)q^{L(0)-c/24}+\sum_{i=1}^\infty
G_{2i+2}(\tau) {\rm tr}|_M o(L[2i]w) q^{L(0)-c/24}.$$ From
(\ref{l20}) and the previous formula applied for
$V=L(-22/5,0)^{\otimes^{12}}$ and $v=v_i$, we get
$${\rm tr}|_M o(v_{i+1}) q^{L^{tot}(0)+11/5}={\rm tr}|_{M} o(L^{tot}[-2]v_{i})q^{L^{tot}(0)+11/5}$$
$$=\Theta \left( {\rm tr}|_M o(v_i) q^{L^{tot}(0)+\frac{11}{5}} \right)+ i(12-i+1)(-\frac{11}{5}G_4(\tau)
){\rm tr}|_M o(v_{i-1}) q^{L^{tot}(0)+11/5},$$ which is equivalent
to the formula (\ref{recdiff}). Furthermore, \bea && \frac{1}{12!}
{\rm tr}|_{M} o(L^{tot}[-2] \cdot
v_{12})q^{L^{tot}(0)+\frac{11}{5}}={\rm tr}|_M o(L[-2]^2{\bf 1}
\otimes  \cdots \otimes L[-2]{\bf
1})q^{L^{tot}(0)+\frac{11}{5}}+\cdots \nn &&+{\rm tr}|_M o(L[-2]
{\bf 1} \otimes \cdots \otimes L[-2]^2{\bf
1})q^{L^{tot}(0)+\frac{11}{5}}=0, \nonumber \eea because of
(\ref{l20}). Here $L^{tot}[-2]$ is a Virasoro generator acting on
the tensor product vertex operator algebra via comultiplication.
Now we can proceed as in Lemma \ref{prelast}, and we get
$R_{12}=0$.

\subsection{On $L(c_{2,5},0)^{\otimes^{12}}$ and $L(c_{2,27},0)$}

In this section we give a combinatorial interpretation of
(\ref{wat2}) in terms of colored partitions and discuss some
related work.

Let us recall \cite{A} that $q^{-11/60}{\rm ch}_1(q)$ (resp.
$q^{1/60}{\rm ch}_2(q)$) is actually the generating series for the
number of partitions in parts congruent to $\pm 2 \ {\rm mod} \
5,$ (resp. $\pm 1 \ {\rm mod \ 5}$). Let
$P_{j_1,j_2,j_3,j_4,j_5}(n)$ denotes the number of colored
partitions of $n$ where every part of size $i \ {\rm mod} \ { 5}$
can be colored in at most $j_i$ colors. Then we have
\begin{proposition} \label{ssss} For every $n \geq 2$,
$$P_{11,1,1,11,0}(n)=11 P_{6,6,6,6,0}(n-1)+P_{1,11,11,1,0}(n-2).$$
\end{proposition}
\noindent {\em Proof:} It is known that the generating functions
of colored partitions in which every part of size $j$ can be
colored with at most $c_j$ colors is given by
$$\prod_{j=1}^\infty \frac{1}{(1-q^j)^{c_j}}.$$
The statement now follows from (\ref{wat2}). \epf

\begin{remark}
{\em Modular forms $\WC$ associated to irreducible characters of
$\mathcal{M}(2,2k+1)$ Virasoro minimal models (essentially
Andrews-Gordon series \cite{FF}) have recently been studied in
\cite{MMO} in connection with supersingular $j$-invariants. We
proved that the quotient $\WC$ is trivial if and only if
$k=6s^2-6s+1$, $s \geq 2$, which is equivalent to a family of
$q$-series identities among irreducible characters. For $s=2$,
($k=13$) the vanishing is equivalent to the following three term
combinatorial identity:
$$P_{27,12}(n)=P_{27,6}(n-1)+P_{27,3}(n-2),$$
where $P_{a,b}(n)$ denotes the number of partitions of $n$ into
parts which are not congruent to $0, \pm b \ {\rm mod} \ a$. This
identity, compared with Proposition \ref{ssss}, indicates that
vertex operator algebra $L(-22/5,0)^{\otimes^{12}}$ shares some
similarities with $L(c_{2,27},0)$. For example, both vertex
operator algebras have exactly $13$ linearly independent
irreducible characters.
}
\end{remark}

\begin{remark} \label{4const}
{\em Four constants $-\frac{11}{5}, \frac{-25}{4}, -15$ and $
-40$, appearing in Lemma \ref{lemmalast} all give rise to
two-dimensional modular invariant spaces coming from irreducible
characters of certain integrable lowest weight representations of
Kac-Moody Lie algebras (e.g., in the $\lambda=-\frac{25}{4}$ case,
for a fundamental system of solutions of $\Theta^2
y-\frac{25}{4}G_4 y=0$ we can take $f_1$ and $f_2$ as in
Proposition \ref{weber}). Similarly for $\lambda=-40$ (see Section
8) and $\lambda=-15$ (cf. \cite{Mu1}). Only the $\lambda=0$ case
has no interpretation in terms of graded dimensions, in which case
for a fundamental system of solutions we can take
$$g_1(\tau)=\int_{\tau}^{i \infty} \eta(s)^4 ds \ \ {\rm and} \  \
g_2(\tau)=1.$$ }
\end{remark}

\section{Supersingular $j$-invariants}

In this section we closely follow \cite{O}. Let $f(\tau) \in M_k$
and $\tilde{F}(f,x)$ as in (\ref{rep}). Also, let
$$h_k(x):= \left\{ \begin{array}{c} \ \ \ \ \ \ \ \ \ \ \ \ \ \ \
1 \ {\rm if} \ k\equiv
0 \pmod{12}, \\  x^2(x-1728) \  {\rm if} \ k\equiv 2 \pmod{12},\\
\ \ \ \ \ \ \ \ \ \ \ \ \ \ \ x \ {\rm if} \ k\equiv 4 \pmod{12}, \\
\ \ \ \ \ x-1728 \ {\rm if} \
k\equiv 6 \pmod{12}, \\ \ \ \ \ \ \ \ \ \ \ \ \ \ \ x^2 \ {\rm if} \ k\equiv 8 \pmod{12}, \\
\  x(x-1728) \ \ {\rm if} \  k\equiv 10 \pmod{12}. \end{array}
\right.
$$
Then, we define the {\it{divisor polynomial}} $F(f,x)$ by
\begin{equation}\label{dp}
F(f,x):=h_k(x)\tilde F(f,x).
\end{equation}
Let us recall a few known results about supersingular
$j$-invariants. We say that an elliptic curve over a field $K$ of
characteristic $p>0$ is {\em supersingular} if the group
$E(\bar{K})$ has no torsion \cite{Si}. It is known that there are
only finitely many supersingular curves over $\bar{\mathbb{F}}_p$.
If $p\geq 5$ is prime, then the supersingular loci $S_p(x)$ and
$\widetilde{S}_p(x)$ are defined in $\mathbb{F}_p[x]$ by the
following products over isomorphism classes of supersingular
elliptic curves:
$$
S_p(x):= \prod_{E/\overline{\mathbb{F}}_p\ {\text {\rm
supersingular}}} (x-j(E)),
$$
\begin{equation}\label{ssplocus}
\widetilde{S}_p(x):= \prod_{\substack{E/\overline{\mathbb{F}}_p \
\ {\text {\rm supersingular}}\\ j(E) \not \in \{0, 1728\}}}
(x-j(E)).
\end{equation}
It is known that the polynomial $S_p(x)$ splits completely in
$\mathbb{F}_{p^2}$ \cite{Si}. Define $\epsilon_{\omega}(p)$ and
$\epsilon_{i}(p)$ by
$$
\epsilon_{\omega}(p):=\begin{cases}
0 &{\text {\rm if}}\ p\equiv 1 \pmod 3,\\
1 &{\text {\rm if}}\ p\equiv 2 \pmod 3,
\end{cases}
$$
$$
\epsilon_{i}(p):=\begin{cases}
0 &{\text {\rm if}}\ p\equiv 1 \pmod 4,\\
1 &{\text {\rm if}}\ p\equiv 3 \pmod 4,
\end{cases}
$$
The following proposition relates $S_p(x)$ to $\tilde S_p(x)$
\cite{Si}.
\begin{proposition}
If $p\ge5$ is prime, then
\begin{displaymath}
\begin{split}
S_p(x)&=x^{\epsilon_{\omega}(p)}(x-1728)^{\epsilon_{i}(p)}\cdot\prod_{\alpha\in \mathfrak{S}_p}(x-\alpha)\cdot\prod_{g\in \mathfrak{M}_p}g(x)\\
&=x^{\epsilon_{\omega}(p)}(x-1728)^{\epsilon_{i}(p)}\tilde S_p(x).
\end{split}
\end{displaymath}
\end{proposition}
\noindent Deligne found the following congruence (see
\cite{Serre}).
\begin{theorem}\label{deltheorem}
If $p\ge5$ is prime, then
$$
F(E_{p-1},x)\equiv S_p(x) \pmod p.
$$
\end{theorem}

\begin{remark} \label{VS}
{\em The Von-Staudt congruences imply for primes $p$, that
$\tfrac{2(p-1)}{B_{p-1}} \equiv 0 \pmod p$, where $B_n$ denotes
the usual $n$th Bernoulli number.  It follows that
\begin{displaymath}
E_{p-1}(\tau)\equiv 1 \pmod p.
\end{displaymath}
If $p\geq 5$ is prime, then Theorem~\ref{deltheorem} combined with
the definition of divisor polynomials, implies that if $f(\tau)\in
M_{p-1}$ and $f(\tau)\equiv 1 \pmod p$, then
\begin{displaymath} \label{SSC}
F(f,j(\tau))\equiv S_p(j(\tau)) \pmod p.
\end{displaymath}
}
\end{remark}

\section{Symmetric powers associated to level one representations of $D_4^{(1)}$}

\noindent In this section we focus on a particular family of
modular forms which give supersingular $j$-invariants in prime
characteristics. In what follows $p \geq 5$ is prime and \be
\label{2m3} p=2m+3. \ee
It is known \cite{K} (see also \cite{Mu1}) that the graded
dimensions of level one highest weight modules for $D_4^{(1)}$
span a two-dimensional vector space $U$ with a basis consisting of
eighth powers of two Weber modular functions:
$$\goth{f}^8=q^{-1/6} \prod_{n=1}^\infty (1+q^{n-1/2})^8$$
$$\goth{f}_2^8=q^{1/3} \prod_{n=1}^\infty (1+q^n)^8.$$
Notice that $\frac{1}{3}-\frac{1}{6}=\frac{1}{6}$, so
$\mathcal{W}_U(\tau)=\eta(\tau)^4$ by \cite{M2}. Thus, (cf.
\cite{M2} or \cite{MMS}):
\begin{lemma} \label{-40}
The infinite products $\goth{f}^8$ and $\goth{f}_2^8$ form a
fundamental set of solutions of the ODE (\ref{ODEl}) with
$\lambda=-40$.
\end{lemma}
We will focus on the $m$-th symmetric power of $U$. As we already
mentioned ${\rm Sym}^m(U)$ is $(m+1)$-dimensional.
In what follows we will use a result  from \cite{KZ}. In that
paper, among other things, Kaneko and Zagier studied the
generating series of the form
$$G_{\alpha}(x)=(1-3E_4 x^4+2E_6 x^6)^\alpha,$$
for some special $\alpha \in \mathbb{Q}$. For $l \in \mathbb{N}$
and $\alpha \in \mathbb{Q}$ let us define
$$G_{l,\alpha}={\rm Coeff}_{x^{2l}}(1-3E_4 x^4+2E_6 x^6)^\alpha \in \mathbb{Q}[E_4,E_6].$$
The following result is from \cite{KZ}.
\begin{proposition} \label{propKZ}
For every prime $p \geq 5$
$$G_{\frac{p-1}{2},\frac{p-3}{6}} \equiv 12^{\frac{p-1}{2}} \
({\rm mod} \ p).$$
\end{proposition}
\noindent The main idea behind the proof of Proposition
\ref{propKZ} is the congruence
\begin{equation} \label{1/2}
(1-3E_4x^4+2E_6 x^6)^{(p-3)/6} \equiv (1-3E_4x^4+2E_6 x^6)^{-1/2} \
({\rm mod} \ p),
\end{equation}
the Von-Staudt congruences (cf. Remark \ref{VS}) and a
parametrization of the elliptic curve $E_\tau$ by using the
Weierstrass $\wp$-function (see \cite{KZ} for details).

Let us recall again that the graded vector space
$M=\mathbb{C}[E_4,E_6]$ admits a graded map $\Theta$
(\ref{thetafirst}) from $M_k$ to $M_{k+2}$, which can be written
as \cite{Z} \be \label{theta} \Theta=-\frac{E_6}{3}
\frac{\partial}{\partial E_4}-\frac{E_4^2}{2}
\frac{\partial}{\partial E_6}. \ee The goal of this section is to
prove the following result.
\begin{theorem} \label{two} Let $V={\rm Sym}^m(U)$. Then
\begin{itemize}
\item[(i)] For every $m \geq 1$,
$${\frac{W'_V(\tau)}{W_V(\tau)}}=(-1)^{m+1} \frac{(m+1)!}{6^{m+1}}
{\rm Coeff}_{x^{2m+2}} (1-3E_4 x^4+2E_6x^6)^{\frac{m}{3}}.$$
\item[(ii)] For $m$ and $p$ as in (\ref{2m3})
$${\frac{W'_V(\tau)}{W_V(\tau)}} \equiv (-1)^{(p-1)/2} \left(
\frac{2}{p} \right)\left(\frac{p-1}{2}\right)! \ ({\rm mod} \
p),$$ where $\left(\frac{\cdot}{p}\right)$ is the
Legendre symbol. \\
\item[(iii)] For $m$ and $p$ as in (\ref{2m3}) $$F(\WC,j(\tau))
\equiv S_p(j(\tau)) \ ({\rm mod} \ p).$$
\end{itemize}
\end{theorem}
\noindent {\em Proof:} Firstly, \be
(1+x+y)^{m/3}=\sum_{r=0,s=0}^\infty
\frac{m/3(m/3-1)\cdots(m/3-r-s+1)}{r! s!} x^r y^s, \ee gives \bea
&&(1-3E_4x^4+2E_6x^6)^{m/3} \nn &&= \sum_{l=0}^\infty \left(
\sum_{r,s \geq 0,2r+3s=l} \frac{m/3(m/3-1)\cdots(m/3-r-s+1)}{r!
s!} (-3E_4)^r (2E_6)^s \right) x^{2l}. \eea Clearly, \bea
G_{l,m/3}&=&{\rm Coeff}_{x^{2l}}(1-3E_4x^4+2E_6x^6)^{m/3} \nn &=&
\sum_{r,s \geq 0, 2r+3s=l} \frac{m/3(m/3-1)\cdots(m/3-r-s+1)}{r!
s!} (-3E_4)^r (2E_6)^s. \eea Now, let
$$\bar{G}_{l,m/3}=\frac{l!}{2^l 3^l} G_{l,m/3}.$$
{\em Claim:} We have $$\bar{G}_{2,m/3}=\frac{-m E_4}{18}, \ \
\bar{G}_{3,m/3}=\frac{m E_6}{54},$$ and for $l=2r+3s \geq 4$, \be
\label{recursionKZ} \bar{G}_{l,m/3}=\Theta
\bar{G}_{l-1,m/3}+(l-1)(m-l+2) \frac{-E_4}{18} \bar{G}_{l-2,m/3}.
\ee

To prove the claim it is enough to consider the coefficient of
$E_4^r E_6^s$ on both sides of (\ref{recursionKZ}) and check the
initial conditions. The coefficient of $E_4^r E_6^s$ on the
left-hand side of (\ref{recursionKZ}) is equal to \be
\label{leftKZ} \frac{(2r+3s)! m/3(m/3-1) \cdots (m/3-r-s+1) (-1)^r
}{r! s! 2^{2r+2s} 2^{r+3s}}. \ee The coefficient of $E_4^rE_6^s$
of the right-hand side of (\ref{recursionKZ}) is  \bea
\label{rightKZ} && \frac{(2r+3s-1)! m/3(m/3-1) \cdots (m/3-r-s+1)
(-1)^r }{r! (s-1)! 2^{2r+3s} 3^{r+3s-1}}\nn && +\frac{(2r+3s-1)!
m/3(m/3-1) \cdots (m/3-r-s+2) (-1)^{r+1} }{(r-2)! s! 2^{2r+2s-1}
3^{r+3s+1}} \nn && +\frac{(2r+3s-2)! m/3 (m/3-1) \cdots
(m/3-r-s+2) (-1)^{r}(2r+3s-1)(m-2r-3s+2) }{(r-1)! s! 2^{2r+2s-1}
3^{r+3s+1}}, \eea where for $r=0$ (resp. $r=1$) the second and
third  (resp. second) term drops. From the identity
\be
\frac{3s(\frac{m}{3}-r-s+1)}{2r+3s}-\frac{2r(r-1)}{3(2r+3s)}+\frac{2r(m-2r-3s+2)}{3(2r+3s)}=\frac{m}{3}-r-s+1,
\ee it follows that (\ref{leftKZ}) is equal to (\ref{rightKZ}).
Thus, the recursion holds. It is easy to see that \be
\label{initial} \bar{G}_{2,m/3}=\frac{-m E_4}{18}=(-40)m {G_4}, \
\ \bar{G}_{3,m/3}=\frac{m E_6}{54}=(-40)m \Theta G_4. \ee Now,
equations (\ref{recursionKZ}) and (\ref{initial}), together with
Lemma \ref{prelast} and (\ref{m1}) imply that $\bar{G}_{l,m/3}$
satisfy the same recursion and the same initial conditions as
$R_{l-1}$ in Lemma \ref{prelast}, for $Q=-40 G_4$. Thus the
formula (i) holds. The part (ii) now follows from Proposition
\ref{propKZ} and $\frac{((p-1)/2)!}{6^{(p-1)/2}}12^{(p-1)/2}
\equiv ((p-1)/2)! \left( \frac{2}{p} \right) \ ({\rm mod} \ p).$
Finally, the equation (iii) follows from (ii) and Remark
\ref{SSC}. \epf

\begin{remark}
{\em Notice that our proof provides also a description of
$G_{l,m/3}$ for {\em every} $l \leq m$ via the recursion in Lemma
\ref{prelast}. Also, for $m \equiv 0 \ ({\rm mod} \ 3)$ the
vanishing of $\W$ is equivalent to (\ref{ex3}). }
\end{remark}

\begin{remark}
{\em It would be nice to have a purely representation theoretic
proof of Theorem \ref{two} via certain differential equations of
order two studied in \cite{KZ} and \cite{KK}, by using techniques
from \cite{M1}. }
\end{remark}

\end{document}